\pdfoutput=1
\documentclass[11pt]{amsart}

\usepackage{amsthm,amssymb,mathtools,comment,colonequals,float,nicematrix,amsmath}

\usepackage[lite,nobysame]{amsrefs}
\usepackage{dsfont}

\usepackage{enumitem}
\setlist[enumerate]{label=$(\alph*)$}

\usepackage{tikz}
\usetikzlibrary{arrows,backgrounds}
\usepackage{tikz-cd}
\tikzstyle{vertex}=[circle, draw, inner sep=0pt, minimum size=6pt]
\newcommand{\vertex}{\node[vertex]}

\newcommand{\R}{{\mathbb R}}
\newcommand{\N}{{\mathbb N}}
\newcommand{\Z}{{\mathbb Z}}
\newcommand{\Q}{{\mathbb Q}}
\newcommand{\ps}{{\mathbb P}}
\newcommand{\A}{{\mathbb A}}
\newcommand{\F}{{\mathbb F}}

\newcommand{\p}{{\mathfrak p}}

\newcommand{\B}{{\mathcal B}}
\newcommand{\J}{{\mathcal J}}

\newcommand\inner[2]{\langle #1, #2 \rangle}
\newcommand{\Mod}[1]{\ \mathrm{mod}\ #1}
	
\DeclareMathOperator\Path{Path}
\DeclareMathOperator\cycle{Cycle}
\DeclareMathOperator\coker{coker}

\DeclareMathOperator\ep{ep}

\DeclareMathOperator\Jac{Jac}
\DeclareMathOperator\height{ht}
\DeclareMathOperator\proj{Proj}
\DeclareMathOperator\spec{Spec}
	
\newcommand\blfootnote[1]{%
	\begingroup
	\renewcommand\thefootnote{}\footnote{#1}%
	\addtocounter{footnote}{-1}%
	\endgroup
}

\numberwithin{equation}{section}
	
\theoremstyle{plain} 
	\newtheorem{theorem}[equation]{Theorem}
	\newtheorem{corollary}[equation]{Corollary}
	\newtheorem{lemma}[equation]{Lemma}
	\newtheorem{proposition}[equation]{Proposition}
	
\theoremstyle{definition}
	\newtheorem{definition}[equation]{Definition}
	
\theoremstyle{remark}
	\newtheorem{remark}[equation]{Remark}
	\newtheorem{example}[equation]{Example}

\begin{document}
		
\title[]{On the proportion of metric matroids whose Jacobians have nontrivial $\mathbf{p}$-torsion}

\author{Sergio Ricardo Zapata Ceballos}
		
\address{}
		
		


\date{}
		
\begin{abstract}
We study the proportion of metric matroids whose Jacobians have nontrivial $p$-torsion. We establish a correspondence between these Jacobians and the $\F_p$-rational points on configuration hypersurfaces, thereby relating their proportions. By counting points over finite fields, we prove that the proportion of these Jacobians is asymptotically equivalent to $1/p$.
\end{abstract}

\maketitle
		

\section{Introduction}
\label{sec:intro}

\blfootnote{MSC2020 \textit{Mathematics Subject Classification}. Primary 11G25, 14G05, 14N20, 60B99; Secondary 05C31, 05C50, 05C76, 14M12,   81Q30.}

\blfootnote{\textit{Key words and phrases}. Arithmetic statistics, regular matroid, series extension, Jacobian group, torsion, configuration polynomial, configuration hypersurface, finite field, rational point, density.}

A graph can be associated with a finite abelian group that plays a significant role in various mathematical disciplines, including arithmetic geometry, combinatorics, and statistical physics. Depending on the context, this group is referred to by different names including: sandpile group, chip-firing group, or critical group. In the context of geometry, it is called Jacobian as it serves as a discrete counterpart to the Jacobian of a Riemann surface or of an algebraic curve over a finite field, highlighting the analogy between Riemann surfaces and graphs~\cite{Bake-Norine07,Baker,Horton}. A more direct connection is found in the group of components of the N\'eron model of a Jacobian of a curve over a local field, which is given as a Jacobian of a graph \cite{Lorenzini89, Bosch02}. We refer the reader to~\cite{Lorenzini08} for a discussion of these and other connections.

These connections, coupled with the importance of the Jacobian as a graph invariant, have driven the exploration of arithmetic statistical problems concerning families of graphs. For instance, Cohen-Lenstra heuristics have been studied for random graphs in \cite{Lorenzini08,Clancy-et15,Wood17}. However, very few families have been studied, and little is known about the structure of the Jacobians. The present work contributes to this area by investigating a new family of graphs.

This paper addresses a variational problem concerning the family of metric graphs. We start with a graph $G$ and assign a positive integer value $\lambda(e)$ to each edge $e$ of $G$. By repeatedly subdividing each edge $\lambda(e)$ times, we obtain a new graph denoted $G_\lambda$, which is referred to as a metric graph. Our first major result (Theorem~\ref{thm:order_of_jacobian_parametrized_by_conf_polynomial}) establishes a link between metric graphs and the configuration polynomial of $G$, also referred to as the graph polynomial or Symanzik polynomial of $G$: 
\begin{equation*}
	\Psi_G(\lambda)=\#\Jac(G_\lambda),
\end{equation*}
where $\Jac(G_\lambda)$ denotes the Jacobian of $G_\lambda$ and $\Psi_G$ is the configuration polynomial of $G$. Motivated by this relationship, we study the proportion of $\Jac(G_\lambda)$ with nontrivial $p$-torsion in the family $\{\Jac(G_\lambda)\}_\lambda$ by counting the zeroes of $\Psi_G$ modulo $p$.

An important consequence of Theorem~\ref{thm:order_of_jacobian_parametrized_by_conf_polynomial} is that the $\F_p$-rational points on the hypersurface $X_G$ cut out by $\Psi_G$ parameterize the $\lambda$'s for which $\Jac(G_\lambda)$ has nontrivial $p$-torsion. This result reveals an arithmetic aspect of $X_G$ and allows us to predict the proportion of Jacobians with non-trivial $p$-torsion by estimating the number of $\F_p$-rational points on $X_G$. Additionally,  Theorem~\ref{thm:order_of_jacobian_parametrized_by_conf_polynomial} provides new insights into $X_G$ contributing to the understanding of its geometry, a subject of ongoing investigation. These hypersurfaces are particularly significant in the study of Feynman integrals in quantum field theory~\cite{Marcolli}. However, to my knowledge, the approach of counting rational points on these hypersurfaces has not been explored.

We develop the theory in the language of matroids as all the aforementioned constructions naturally extend to regular matroids. After introducing regular matroids and their Jacobians, we generalize the construction of the family of metric graphs to matroids through systematic serial extensions and direct sums, referring to these as metric matroids. We establish several general properties of metric matroids, the most pertinent being that a metric matroid associated with a regular matroid is itself regular (Proposition~\ref{prop:M_lambda_is_regular}).

Next, we introduce the configuration hypersurface $X_M$ of a regular matroid $M$ and explore some of its geometric properties. We establish bounds for $\#X_M(\F_p)$, enhancing the existing bounds in the literature for this particular case (Theorem~\ref{thm:rational_points_big_O}). Finally, we leverage these results to prove Theorems~\ref{thm:limit_of_height_of_lambdas} and \ref{thm:estimate_of_limit}, showing that the proportion of Jacobians with nontrivial $p$-torsion in the family of metric matroids is asymptotically equivalent to $1/p$.


\section{background}
\label{sec:background}

\subsection{Regular matroids} 
\label{subsec:RegMtrds}
We assume that the reader is familiar with the basic theory of matroids; a standard reference is the book on matroids by Oxley~\cite{Oxley11}.

A matrix $A$ over the integers is said to be \emph{totally unimodular} if the determinant of any square submatrix of $A$ is in $\{-1,0,1\}$.

A matroid $M$ is called \emph{regular} if it can be represented over $\Q$ by an $r\times n$ totally unimodular matrix of rank $r$. 

In what follows, by a graph we will mean a finite undirected multi-graph. One can associate a regular matroid $M(G)$ to a graph $G$, by taking as a totally unimodular representation the matrix resulting from removing some appropriate rows of the incidence matrix of $G$. A matroid $M$ is called \emph{graphic} if $M=M(G)$ for some graph $G$.

A matroid $M$ is said to be \emph{irreducible} (or \emph{connected}) if it cannot be written as the direct sum of two nonempty matroids. In particular, it is well-known that if $G$ is a loopless graph without isolated vertices and $\#V(G)\geq 3$, then $M(G)$ is irreducible if and only if $G$ is 2-connected.

If $M$ and $N$ are two matroids, then we say that $N$ is a \emph{series extension} of $M$ if there is a cocircuit $\{e,f\}$ of $N$ such that $e\neq f$ and $N/e =M$. In particular, one has that series extensions of regular matroids are regular as well (see~\cite{Zapata}*{Proposition $1.5.4$}).

\subsection{Jacobian groups}
\label{subsec:JacGrps}
Let $M$ be a regular matroid on $E$ represented over $\Q$ by an $r\times n$ totally unimodular matrix $A$ of rank $r$. Let us consider the $\Q$-vector space $\Q^n$ equipped with the canonical inner product. We define $\Lambda_A(M)\colonequals\ker A\cap \Z^n$. This is a full rank integral $\Z$-lattice of $\ker A$ so that we can speak of its \emph{dual} $\Lambda_A(M)^\#$, that is
\begin{equation*}
	\Lambda_A(M)^\#=\{x\in\ker A : \inner{x}{y}\in\Z~\forall~y\in\Lambda_A(M)\}.
\end{equation*}
It turns out that the isometry classes of these lattices are independent of the totally unimodular representation of $M$ (see~\cite{Merino}*{$\S 4.2	$}). The \emph{Jacobian} $\Jac(M)$ is defined to be the determinant group of the lattice $\Lambda_A(M)$, i.e., $\Jac(M)=\Lambda_A(M)^\#/\Lambda_A(M)$.	The order of $\Jac(M)$ equals the number of bases of $M$ (cf.~\cite{Merino}*{Theorem $4.3.2$}). In fact, as a consequence of Lemma $1$ from Section $4$ in~\cite{Bacher-et97}, there is a natural isomorphism between $\Jac(M)$ and $\coker(AA^t)$ (cf.~\cite{Zapata}*{$\S 2.1.3$}). In particular, $\#\Jac(M)=\det(AA^t)$.

\section{Metric matroids}
\label{sec:metric_matroids}
In the sequel, the set of positive integers will be denoted by $\N$. We first review the notion of metric graphs; we then give our generalization to matroids. 

\subsection{Metric graphs}    	
\label{subsec:metric_graphs}	
Let $G$ be a graph with edge set $E(G)$. Consider the following operation on $G$. Given an arbitrary map $\lambda\colon E(G)\to \N$, we let $G_\lambda$ be the graph obtained from $G$ by replacing every  edge $e\in E(G)$ with a path of length $\lambda(e)$. The graph $G_\lambda$ is called a \emph{metric graph} and the pair $(G,\lambda)$ is called a \emph{model} of $G_\lambda$. Intuitively, the map $\lambda$ assigns a length to each edge of $G$.

Through this process, we derive a family of graphs $\{G_\lambda\}_{\lambda\in \N^{E(G)}}$ from the original graph $G$. It is noteworthy that the original graph $G$ itself is a member of this family, represented as $G_\mathbf{1}$, where $\mathbf{1}$ denotes the constant function $1$.

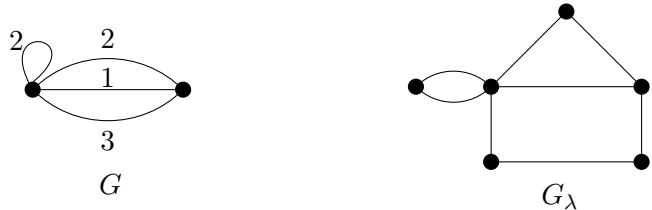
\begin{figure}[H]
	\begin{minipage}{0.5\textwidth}
		\centering
		\begin{tikzpicture}
			\draw[fill] (-1,0) circle[radius=.1];
			\draw[fill] (1,0) circle[radius=.1];
			\draw[every loop/.style={distance=10mm}]  (-1,0) to[in=45,out=130,loop,rotate=-6] node[left=2]{$2$} (1,2);
			\draw(-1,0) to [bend left=45] node[above]{$2$}(1,0);
			\draw(-1,0) to node[above=-3]{$1$}(1,0);
			\draw(-1,0) to [bend right=45] node[below]{$3$}(1,0);
		\end{tikzpicture}
		
		\hspace{15pt}$G$
	\end{minipage}
	\begin{minipage}{0.4\textwidth}
		\centering
		\begin{tikzpicture}
			\draw[fill] (-2,0) circle[radius=.1];
			\draw[fill] (-1,0) circle[radius=.1];
			\draw[fill] (1,0) circle[radius=.1];
			\draw[fill] (0,1) circle[radius=.1];
			\draw[fill] (-1,-1) circle[radius=.1];
			\draw[fill] (1,-1) circle[radius=.1];
			\draw(-2,0) to [bend right=45] (-1,0);
			\draw(-2,0) to [bend left=45] (-1,0);
			\draw(-1,0) to  (1,0);
			\draw(-1,0) to  (-1,-1);
			\draw(-1,-1) to  (1,-1);
			\draw(1,0) to  (1,-1);
			\draw(-1,0) to  (0,1);
			\draw(1,0) to  (0,1);
		\end{tikzpicture}
		
		\hspace{20pt} $G_\lambda$
	\end{minipage}
	\caption{Example of a metric graph.}	
\end{figure}

\begin{example}
	Let $n\in\N$. The \emph{$n$-path graph}, denoted by $\Path_n$, is the graph with vertex set  $V=\{0,\dots,n\}$, edge set $E=\{e_1,\dots,e_n\}$, and endpoint map given by $\ep(e_i)=\{i-1,i\}$ for all $i\in\{1,\dots,n\}$.
	If $G=\Path_1$ and $\mathbf{n}$ is the map that sends $e_1$ to $n$, then $G_\mathbf{n}=\Path_n$ for all $n\in\N$.
\end{example}

\begin{example}
	Let $n\in\N$. The \emph{$n$-cycle graph}, denoted by $\cycle_n$, is the graph with vertex set  $V=\Z/n\Z$, edge set $E=\{e_1,\dots,e_n\}$, and endpoint map defined by the rule $\ep(e_{i+1})=\{\overline{i},\overline{i+1}\}$ for all $\overline{i}\in\Z/n\Z$. If $G=\cycle_1$, then $G_\mathbf{n}=\cycle_n$ for all $n\in\N$.
\end{example}
		
\subsection{Metric matroids}
\label{subsec:metric_matroids}			
Consider a matroid $M$ on $E$ of rank $r$ with $n$ elements and fix a map $\lambda\in\N^E$. For each $e\in E$, we denote by $E_e(M)$ a set containing $\lambda(e)$ elements such that $e\in E_e(M)$ and $E_e(M)\cap E_f(M)=\emptyset$ for $e\neq f$. We define
\begin{equation*}
	E(M_\lambda)\colonequals\bigsqcup_{\mathclap{e\in E(M)}}E_e(M).
\end{equation*}
Let $\B(M_\lambda)$ be the collection of all subsets $B$ of $E(M_\lambda)$ for which there exists a basis $B'$ of $M$ and a tuple $(x_e)_{e\notin B'}\in \prod_{e \notin B'}E_e(M)$ such that $B=E(M_\lambda)-\{x_e\vcentcolon e\notin B'\}$.

Equivalently, there exists a basis $B''$ of $M^\ast$ and a tuple $(x_e)_{e\in B''}\in \prod_{e \in B''}E_e(M)$ such that $B=E(M_\lambda)-\{x_e\vcentcolon e\in B''\}$, where $M^\ast$ is the dual matroid of $M$.

\begin{proposition}\label{prop:M_lambda_is_a_matroid}
	The collection $\B(M_\lambda)$ is the set of bases of a matroid on $E(M_\lambda)$.
\end{proposition}

We need the following fact to prove Proposition~\ref{prop:M_lambda_is_a_matroid}.

\begin{lemma}\label{lem:property_B2_star}
	Let $\B$ be the set of bases of a matroid $M$. If $B_1,B_2\in\B$ and $x\in B_2-B_1$, then there is $y\in B_1-B_2$ such that $(B_1-\{y\})\cup\{x\}\in\B$.
\end{lemma}
\begin{proof}
	See Lemma $2.1.2$ in \cite{Oxley11}.
\end{proof}

\begin{proof}[Proof of Proposition~\ref{prop:M_lambda_is_a_matroid}]
	We must prove that $\B(M_\lambda)$ satisfies the following two conditions:
	\begin{enumerate}
		\item $\B(M_\lambda)\neq\emptyset$; 
		
		\item if $B_1,B_2\in \B(M_\lambda)$ and $x\in B_1-B_2$, then there is $y\in B_2-B_1$ such that $(B_1-\{x\})\cup\{y\}\in\B(M_\lambda)$.
	\end{enumerate}
	The first condition follows by definition as the collection of bases of $M^\ast$ is nonempty. 
	
	To prove $(b)$, we pick $B_1, B_2\in \B(M_\lambda)$ and suppose $x\in B_1-B_2$. By definition, there are bases $B'_1=\{a_1,\dots,a_{n-r}\}$ and $B'_2=\{b_1,\dots,b_{n-r}\}$ of $M^\ast$  and tuples $(c_1,\dots,c_{n-r})\in \prod_{i=1}^{n-r}E_{a_i}(M)$ and $(d_1,\dots,d_{n-r})\in \prod_{i=1}^{n-r}E_{b_i}(M)$ such that  $B_1=E(M_\lambda)-\{c_1,\dots,c_{n-r}\}$ and $B_2=E(M_\lambda)-\{d_1,\dots,d_{n-r}\}$.
	
	Observe that $B_1-B_2=\{d_1,\dots,d_{n-r}\}-\{c_1,\dots,c_{n-r}\}$ and  $B_2-B_1=\{c_1,\dots,c_{n-r}\}-\{d_1,\dots,d_{n-r}\}$. Thus, there is $j\in\{1,\dots,n-r\}$ such that $x=d_j$. Let us consider two cases. 
	
	\textit{Case 1}. Suppose there exists $k\in\{1,\dots,n-r\}$ such that $b_j=a_k$. In this case, $c_k\in E_{b_j}(M)$ so that $c_k\in B_2-B_1$ and
	\begin{equation*}
		(B_1-\{d_j\})\cup\{c_k\}=E(M_\lambda)-(\{c_1,\dots,c_{n-r},d_j\}-\{c_k\})
	\end{equation*}
	  where
	\begin{equation*}
		(c_1,\dots,c_{k-1},d_j,c_{k+1},\dots,c_{n-r})\in\prod_{i=1}^{n-r}E_{a_i}(M).
	\end{equation*}
	Hence $(B_1-\{d_j\})\cup\{c_k\}\in\B(M_\lambda)$.
	
	\textit{Case 2}. Suppose $b_j\neq a_i$ for all $i\in\{1,\dots,n-r\}$. Thus, we have $b_j\in B'_2-B'_1$.  By Lemma~\ref{lem:property_B2_star}, there exists $a_k\in B'_1-B'_2$ such that $(B'_1-\{a_k\})\cup\{b_j\}$ is a basis of $M^\ast$. Since $a_k\notin B'_2$, we have $E_{a_k}(M)\neq E_{b_i}(M)$ for all $i\in\{1,\dots,n-r\}$. Therefore $c_k\notin E_{b_i}(M)$ for all $i\in\{1,\dots,n-r\}$ and consequently $c_k\in B_2-B_1$. Observe that
	\begin{equation*}
		(B_1-\{d_j\})\cup\{c_k\}=E(M_\lambda)-(\{c_1,\dots,c_{n-r},d_j\}-\{c_k\})
	\end{equation*}
	  where
	\begin{equation*}
		(d_j,c_1,\dots,c_{k-1},c_{k+1},\dots,c_{n-r})\in E_{b_j}(M)\times\prod_{\substack{{i=1}\\i\neq k}}^{n-r}E_{a_i}(M).
	\end{equation*}
	Then $(B_1-\{d_j\})\cup\{c_k\}\in\B(M_\lambda)$ since $\{a_1,\dots,a_{n-r},b_j\}-\{a_k\}$ is a basis of $M^\ast$. This completes the proof.
\end{proof}

\begin{definition}\label{def:M_lambda}
	Let $M$ be a matroid on $E$. If $\lambda\in\N^E$, then the matroid  $(E(M_\lambda),\B(M_\lambda))$, denoted by $M_\lambda$, is called a \emph{metric matroid}.
\end{definition}

The reason $M_\lambda$ is called a metric matroid is that $M(G)_\lambda=M(G_\lambda)$ for any graph $G$ (see Proposition~\ref{prop:M_lambda_for_graphic_matroids}).

\begin{example}
The matroid $M_\mathbf{1}$ given by the constant function $e\mapsto 1$ is equal to $M$. 
\end{example}

\begin{example}
Let $M$ be a matroid on $E$. For $e\in E$, let $\chi_e\colon E\to \{0,1\}$ be the characteristic function of $e$.	The matroid $M_{\mathbf{1}+\chi_e}$ has ground set	 $E(M_{\mathbf{1}+\chi_e})=E\sqcup\{e'\}$ and its collection of bases is
	\begin{equation}\label{eq:bases_of_M_chi_e}
		\B(M_{\mathbf{1}+\chi_e})=\{B\sqcup\{e'\}\vcentcolon B\in\B(M)\}\cup\{B\cup\{e\}\vcentcolon B\in\B(M), e\notin B\},
	\end{equation}
where $\B(M)$ is the collection of bases of $M$.
\end{example}

Let $m$ and $n$ be two nonnegative integers with $m\leq n$. The \emph{uniform matroid} of rank $m$ on a set of $n$ elements, denoted $U_{m,n}$, is the matroid with ground set $E(U_{m,n})$ a set of size $n$ and collection of bases those subsets of $E(U_{m,n})$ whose cardinality is $m$.

\begin{proposition}\label{prop:M_lambda_is_series_extension}
	Let $M$ be a matroid on $E$ and let $e\in E$. Then 
	\begin{enumerate}
		\item\label{prop-a} if $e$ is not a coloop of $M$, then $M_{\mathbf{1}+\chi_e}$ is a series extension of $M$;
		
		\item\label{prop-b} if $e$ is a coloop of $M$, then $M_{\mathbf{1}+\chi_e}\cong M\oplus U_{1,1}$. 
	\end{enumerate}
\end{proposition}
\begin{proof}
	\ref{prop-a} If $e$ is not a coloop of $M$, then from the description of $\B(M_{\mathbf{1}+\chi_e})$ in \eqref{eq:bases_of_M_chi_e}, it follows that the set $\{e,e'\}$ is a cocircuit of $M_{\mathbf{1}+\chi_e}$ satisfying $e\neq e'$. It also follows from~\eqref{eq:bases_of_M_chi_e} that $M_{\mathbf{1}+\chi_e}/e'=M$. This completes the proof.
	
	\ref{prop-b} If $e$ is a coloop of $M$, then
	\begin{equation*}
		M_{\mathbf{1}+\chi_e}=\Big(E\sqcup\{e'\},\{B\sqcup\{e'\}\vcentcolon B\in\B(M)\}\Big),
	\end{equation*}
	which is isomorphic to $M\oplus U_{1,1}$.
\end{proof}

Proposition~\ref{prop:M_lambda_is_series_extension} shows how $M_\lambda$ can be obtained by successively constructing series extensions starting with $M$ and/or adding copies of $U_{1,1}$. If $M$ is irreducible, then $M$ has no coloops; therefore, $M_\lambda$ is obtained from $M$ through series extensions only.

\begin{proposition}\label{prop:rank_of_M_lambda}
	Let $M$ be a matroid on $E$. For any $\lambda\in\N^E$, the rank of $M_\lambda$ is
	\begin{equation*}
		r(M_\lambda)=r(M)+\sum_{\mathclap{e\in E}}(\lambda(e)-1).
	\end{equation*}
\end{proposition}
\begin{proof}
	It is well-known that  $r(M_\lambda)+r(M^\ast_\lambda)=\#E(M_\lambda)$. By the definition of $\B(M_\lambda)$, we have $r(M_\lambda^\ast)=r(M^\ast)$, then $r(M^\ast_\lambda)=\#E-r(M)$. Hence
	\begin{equation*}
		r(M_\lambda)=\#E(M_\lambda)-\#E+r(M)=\sum_{e\in E}(\lambda(e)-1)+r(M).
	\end{equation*}
\end{proof}

\begin{proposition}\label{prop:number_of_bases_of_M_lambda}
	Let $M$ be a matroid on $E$. If $\lambda\in\N^E$, then the number of bases of $M_\lambda$ is
	\begin{equation*}
		\sum_{B\in\B(M)}\prod_{e\notin B}\lambda(e).
	\end{equation*}
\end{proposition}
\begin{proof}
	By definition, there is a bijection between $\B(M_\lambda)$ and the set
	\begin{equation*}
		\bigsqcup_{B\in\B(M)}\prod_{e\notin B}E_e(M).
	\end{equation*}
	Hence the result follows.
\end{proof}

The next proposition will allow us to prove properties about $M_\lambda$ by induction on $\lambda$.

\begin{proposition}\label{prop:inductive_step}
	Let $M$ be a matroid on $E$. Suppose $\lambda\in\N^E$ is such that $\lambda(e)>1$ for some $e\in E$. Then $	(M_{\lambda-\chi_e})_{\mathbf{1}+\chi_e}=M_\lambda$.
\end{proposition}
\begin{proof}
	Suppose that $E(M_\lambda)=\bigsqcup_{f\in E}E_f(M)$. By definition, the ground set of $M_{\lambda-\chi_e}$ can be taken to be  $E(M_{\lambda-\chi_e})=\bigsqcup_{f\in E}E'_f(M)$, where $E'_f(M)=E_f(M)$ for $f\neq e$ and $E'_e(M)=E_e(M)-\{e'\}$ for some $e'\in E_e(M)-\{e\}$. By \eqref{eq:bases_of_M_chi_e}, we have  $E((M_{\lambda-\chi_e})_{\mathbf{1}+\chi_e})=E(M_{\lambda-\chi_e})\sqcup\{e'\}$ and $\B((M_{\lambda-\chi_e})_{\mathbf{1}+\chi_e})$ is the collection
	\begin{equation}\label{eq:proof}
		\{B\cup\{e'\}\vcentcolon B\in\B(M_{\lambda-\chi_e})\}\cup\{B\cup\{e\}\vcentcolon B\in\B(M_{\lambda-\chi_e}), e\notin B\}.
	\end{equation}
	We will show that $\B((M_{\lambda-\chi_e})_{\mathbf{1}+\chi_e})\subseteq \B(M_\lambda)$ and $\B(M_\lambda)\subseteq \B((M_{\lambda-\chi_e})_{\mathbf{1}+\chi_e})$. To see the first containment, we pick an arbitrary basis $B'\in \B((M_{\lambda-\chi_e})_{\mathbf{1}+\chi_e})$. According to \eqref{eq:proof}, we will consider two cases for $B'$.

	\textit{Case 1.}  If $B'=B\cup\{e'\}$ with $B\in\B(M_{\lambda-\chi_e})$, then there exist $B_0\in \B(M)$ and $(x_f)_{f\notin B_0}\in \prod_{f \notin B_0}E'_f(M)$ such that $B=E(M_{\lambda-\chi_e})-\{x_f\vcentcolon f\notin B_0\}$. Since $e'\notin\{x_f\vcentcolon f\notin B_0\}$ we have
	\begin{equation*}
			B'=B\cup\{e'\}=E(M_\lambda)-\{x_f\vcentcolon f\notin B_0\}\in\B(M_\lambda).
	\end{equation*}
		
	\textit{Case 2.} If $B'=B\cup\{e\}$ with $B\in\B(M_{\lambda-\chi_e})$ and $e\notin B$, then there exist a basis $B_0\in \B(M)$ and a tuple $(x_f)_{f\notin B_0}\in \prod_{f \notin B_0}E'_f(M)$ such that $B=E(M_{\lambda-\chi_e})-\{x_f\vcentcolon f\notin B_0\}$. Since $e\notin B$, then $x_e=e$ and $e\notin B_0$. For $f\notin B_0$ with $f\neq e$ define $y_f=x_f$ and $y_e=e'$. Then $(y_f)_{f\notin B_0}\in\prod_{f \notin B_0}E_f(M)$ and
	\begin{equation*}
			B'=B\cup\{e\}=E(M_\lambda)-\{y_f\vcentcolon f\notin B_0\}\in \B(M_\lambda).
	\end{equation*}
	This shows that $\B( (M_{\lambda-\chi_e})_{\mathbf{1}+\chi_e})\subseteq \B(M_\lambda)$.
	
	Now pick $B\in\B(M_\lambda)$. There exist $B_0\in \B(M)$ and $(x_f)_{f\notin B_0}\in \prod_{f \notin B_0}E_f(M)$ such that $B=E(M_\lambda)-\{x_f\vcentcolon f\notin B_0\}$. If $e'\in B$, then $e'\notin \{x_f\vcentcolon f\notin B_0\}$; thus
	\begin{align*}
		B&=E(M_\lambda)-\{x_f\vcentcolon f\notin B_0\}\\
		&=\left( E(M_{\lambda-\chi_e})-\{x_f\vcentcolon f\notin B_0\}\right) \cup\{e'\}\in \B\left( (M_{\lambda-\chi_e})_{\mathbf{1}+\chi_e}\right).
	\end{align*}
	If $e'\notin B$, then $e'\in \{x_f\vcentcolon f\notin B_0\}$, that is, $e\notin B_0$ and $x_e=e'$. We define $y_f=x_f$ if $f\neq e$ and $y_e=e$. So $(y_f)_{f\notin B_0}\in\prod_{f \notin B_0}E'_f(M)$ and 
	\begin{align*}
		B&=E(M_\lambda)-\{x_f\vcentcolon f\notin B_0\}\\
		&=\left( E(M_{\lambda-\chi_e})-\{y_f\vcentcolon f\notin B_0\}\right) \cup\{e\}\in \B\left( (M_{\lambda-\chi_e})_{\mathbf{1}+\chi_e}\right).
	\end{align*}
	Thus $B(M_\lambda)\subseteq \B( (M_{\lambda-\chi_e})_{\mathbf{1}+\chi_e})$.
\end{proof}

Suppose $M$ is a matroid on $E=\{e_1,\dots,e_n\}$. We identify $\lambda\in\N^E$ with the tuple $(\lambda(e_1),\dots,\lambda(e_n))\in\N^n$, and we order $\N^n$ by the lexicographical order.

\begin{proposition}\label{prop:M_lambda_is_regular}
	If $M$ is irreducible (resp. linear or regular), so is $M_\lambda$ for any $\lambda\in\N^E$.
\end{proposition}
\begin{proof}
	We use induction on $\lambda$. The base case is $\lambda=(1,1,\dots,1)\in\N^n$, then $M_\lambda=M$. Thus, the result follows by assumption. We now suppose $\lambda>(1,\dots,1)$ and  that the proposition is true for any $\gamma <\lambda$. 
	
	Let $i=\min\{j\colon\lambda(e_j)>1\}$. Set $\gamma=\lambda-\chi_{e_i}$ so that $\gamma <\lambda$. By the induction hypothesis, the matroid $M_\gamma$ is irreducible (resp. linear or regular). By Proposition~\ref{prop:inductive_step}, $(M_\gamma)_{1+\chi_{e_i}}=M_\lambda$. 
	
	If $M_\gamma$ is irreducible, then $e_i$ is not a coloop of $M_\gamma$, therefore, by Proposition~\ref{prop:M_lambda_is_series_extension}, $M_\lambda$ is a series extension of $M_\gamma$. As a series extension of an irreducible matroid is irreducible (see~\cite{Zapata}*{Proposition $1.4.7$}), the result follows.
	
	If $M_\gamma$ is linear (resp. regular), then we consider two cases.    

	\textit{Case 1.} If $e_i$ is not a coloop of $M_\gamma$, then by Proposition~\ref{prop:M_lambda_is_series_extension}, the matroid $M_\lambda$ is a series extension of $M_\gamma$. As a series extension of a linear (resp. regular) matroid is linear (resp. regular) (see~\cite{Zapata}*{Proposition $1.5.4$}), the result is clear.
		
	\textit{Case 2.} If $e_i$ is a coloop of $M_\gamma$, then $M_\lambda\cong M_\gamma\oplus U_{1,1}$ (see Proposition~\ref{prop:M_lambda_is_series_extension}). Since the direct sum of linear (resp. regular) matroids is linear (resp. regular), the result is clear. 
	  
\end{proof}

\begin{proposition}\label{prop:M_lambda_for_graphic_matroids}
	If $G$ is a graph and  $\lambda\in \N^{E(G)}$, then $M(G)_\lambda=M(G_\lambda)$.
\end{proposition}
\begin{proof}
	First of all, note that $M(G)_{\mathbf{1}+\chi_e}=M(G_{\mathbf{1}+\chi_e})$ for all $e\in E(G)$ (this follows from \eqref{eq:bases_of_M_chi_e}).
	
	We proceed by induction on $\lambda$. If $\lambda=(1,\dots,1)$, then there is nothing to show. Suppose that $(1,\dots,1)<\lambda$ and that the proposition is true for all $\gamma<\lambda$. Let $i=\min\{j\colon \lambda(e_j)> 1\}$ (we are assuming $E(G)=\{e_1,\dots,e_n\}$). Set $\gamma=\lambda-\chi_{e_i}$ so that $\gamma<\lambda$. Thus, by the induction hypothesis, $M(G)_\gamma=M(G_\gamma)$. As a result,
	\begin{equation*}
		M(G)_\lambda=(M(G)_\gamma)_{\mathbf{1}+\chi_{e_i}}=M(G_\gamma)_{\mathbf{1}+\chi_{e_i}}=M((G_\gamma)_{\mathbf{1}+\chi_e})=M(G_\lambda).
	\end{equation*}
\end{proof}


\section{Configuration polynomials}		
\label{sec:config_poly}			
Throughout this section, let $M$ denote a regular matroid on $E$ and $\B(M)$ the collection of bases of $M$. Consider the set of variables $\{\lambda_e\vcentcolon e\in E\}$  indexed by the elements of $M$. The \emph{configuration polynomial} of $M$ is
\begin{align}\label{eq:def_conf_polynomial}
	\Psi_M\colonequals\sum_{\mathclap{B\in\B(M)}}~\prod_{~e \notin B}\lambda_e\in \Z[\lambda_e\vcentcolon e\in E].
\end{align}
By convention, a product of variables indexed by the empty set equals $1$. If $G$ is a graph, we define the configuration polynomial of $G$ to be $\Psi_G\colonequals\Psi_{M(G)}$.

The configuration polynomial of a linear matroid can be defined by the theory of configurations of vector spaces; however, this polynomial depends on the choice of the field of definition. Under this new definition, the configuration polynomial of a regular matroid is independent of the field of definition and has exactly the form as in~\eqref{eq:def_conf_polynomial}. For further discussion of the general case we refer the reader to~\cite{Denham-et21}.

\begin{remark}\label{rmk:general_properties_config_polynomial}
	\
	\begin{enumerate}
		\item By definition, $\Psi_M$ is a homogeneous polynomial of degree $r(M^\ast)$. It is linear in each variable and its coefficients are all $1$.
		
		\item If $e$ is a loop of $M$, then $\lambda_e$ divides $\Psi_M$. Thus, the configuration polynomial of a regular matroid is not irreducible in general.
		
		\item If $e$ is a coloop of $M$, then $\partial_{\lambda_e}\Psi_M=0$. 
	\end{enumerate}
\end{remark}

\begin{example}\label{ex:config_poly_diamond_graph}
	Consider the diamond graph
	\begin{figure}[ht]
		\centering
		\begin{tikzpicture}
			\draw[fill] (1,0) circle[radius=.1];
			\draw[fill] (-1,0) circle[radius=.1];
			\draw[fill] (0,1) circle[radius=.1];
			\draw[fill] (0,-1) circle[radius=.1];
			\draw(1,0) to node[above] {$e_5$} (-1,0);
			\draw(1,0) to node[right] {$e_1$}(0,1);
			\draw(1,0) to node[right] {$e_3$}(0,-1);
			\draw(-1,0) to node[left] {$e_2$}(0,1);
			\draw(-1,0) to node[left] {$e_4$}(0,-1);
		\end{tikzpicture}
		\caption{Diamond graph.}
		\label{fig:diamond_graph}
	\end{figure}
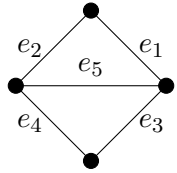
	
	Its configuration polynomial is
	\begin{equation*}
		\lambda_{e_1}\lambda_{e_3}+\lambda_{e_1}\lambda_{e_4}+\lambda_{e_2}\lambda_{e_3}+\lambda_{e_2}\lambda_{e_4}+\lambda_{e_1}\lambda_{e_5}+\lambda_{e_2}\lambda_{e_5}+\lambda_{e_3}\lambda_{e_5}+\lambda_{e_4}\lambda_{e_5}.
	\end{equation*}
\end{example}

\begin{example}\label{ex:config_poly_cycle_graph}
	Let $G=\cycle_n$. Then
	\begin{equation*}
		\Psi_G=\sum_{i=1}^{n}\lambda_{e_i}.
	\end{equation*}
\end{example}

The following proposition relates the irreducibility of a matroid and the irreducibility of its  configuration polynomial. This result was proved in \cite{Denham-et21}. 

\begin{proposition}\label{prop:irreducibility_conf_poluynomial}
	Let $M$ be a regular matroid on $E$.
	\begin{enumerate}
		\item\label{irred-a} If $M_1,\dots,M_n$ are the irreducible components of $M$, then $\Psi_M=\Psi_{M_1}\cdots\Psi_{M_n}$.
		
		\item\label{irred-b} If $M$ is nonempty and has no coloops, then $M$ is irreducible if and only if $\Psi_M$ is irreducible over any field.
		
		\item\label{irred-c} $\Psi_M=1$ if and only if $r(M)=\#E$ if and only if $M$ is isomorphic to a finite direct sum of copies of $U_{1,1}$. 	 	
	\end{enumerate}
\end{proposition}
\begin{proof}
   For proofs of~\ref{irred-a} and \ref{irred-b} see~\cite{Denham-et21}*{Proposition $3.8$}.
   
    \ref{irred-c} If $\Psi_M=1$, then $E-B=\emptyset$ for any $B\in \B(M)$. Therefore, $r(M)=\#E$. If $r(M)=\#E$, then $E$ is the only basis of $M$. Then $M=\bigoplus_{e\in E}M|\{e\}$, where $M|\{e\}$ is the restriction of $M$ to $e$. Since $M|\{e\}\cong U_{1,1}$, $M$ is isomorphic to a finite direct sum of copies of $U_{1,1}$. (We assume that a direct sum of matroids over the empty set equals the empty matroid). Lastly, if $M$ is isomorphic to a finite direct sum of copies of $U_{1,1}$, then $E$ is the only basis of $M$, hence $\Psi_M=1$.
\end{proof}

\begin{proposition}\label{prop:deletion_contraction_conf_polynomial}
	If $M$ is a regular matroid on $E$ and $e\in E$, then
	\begin{equation*}
		\Psi_M=\begin{dcases*}
			\Psi_{M/ e}, & if $e$ is a coloop;\\
			\lambda_e\Psi_{M\backslash e}, & if $e$ is a loop;\\
			\lambda_e\Psi_{M\backslash e}+\Psi_{M/e}, & otherwise.
		\end{dcases*}
	\end{equation*}
\end{proposition}
\begin{proof}
	Cf.~\cite{Denham-et21}*{Proposition $3.12$}.
\end{proof}

\begin{remark}
	When $M$ is irreducible, it has neither loops nor coloops. Therefore, $\Psi_M=\lambda_e\Psi_{M\backslash e}+\Psi_{M/e}$ for any $e\in E$.
\end{remark}

\begin{theorem}\label{thm:order_of_jacobian_parametrized_by_conf_polynomial}
	If $M$ is a regular matroid on $E$, then $\Psi_M(\lambda)=\#\Jac(M_\lambda)$ for all $\lambda\in \N^E$.
\end{theorem}
\begin{proof}
	If $\lambda\in\N^E$, then
	\begin{equation*}
		\Psi_M(\lambda)=\sum_{\mathclap{B\in\B(M)}}~\prod_{~e \notin B}\lambda(e).
	\end{equation*}
	By Proposition~\ref{prop:number_of_bases_of_M_lambda}, $\Psi_M(\lambda)$ equals the number of bases of $M_\lambda$. Since $M_\lambda$ is regular (see Proposition~\ref{prop:M_lambda_is_regular}), we can speak of its Jacobian, whose order is precisely the number of bases of $M_\lambda$. Thus, $\Psi_M(\lambda)=\#\Jac(M_\lambda)$.
\end{proof}

\section{Configuration hypersurfaces}\label{sec:conf_hypersurfaces}

\begin{definition}
	Let $M$ be a regular matroid on $E$ and $k$ be an arbitrary field. The scheme
	\begin{equation*}
		X_M\colonequals\proj k[\lambda_e\vcentcolon e\in E]/(\Psi_M)\subseteq \ps^{\#E-1}_k
	\end{equation*}	
	is called the \emph{configuration hypersurface} of $M$. 
\end{definition}

\noindent If $G$ is a graph, then we define the configuration hypersurface of $G$ to be $X_G\colonequals X_{M(G)}$.

Computations of certain Feynman integrals by Broadhurst and Kreimer~\cite{Kreimer} revealed the presence of multiple zeta values that can be interpreted as periods of configuration hypersurfaces arising from graphs. This connection has drawn significant attention from algebraic geometers, motivating the study of the geometry of these hypersurfaces from a cohomological perspective. Kontsevich conjectured that these hypersurfaces are mixed Tate. Although this conjecture was later disproved by Belkale and Brosnan~\cite{Belkale}, it has underscored the complexity of the singularities associated with configuration polynomials, sparking greater interest in their study. For instance, it has been shown that their singular loci have low codimension~\cite{Denham-et21, Patterson}. For a comprehensive reading on this topic, see~\cite{Marcolli}.

In this section, we will study the geometry of hypersurfaces cut out by homogeneous polynomials that are linear in one of their variables. Note that the configuration polynomial of a regular matroid $M$ is a homogeneous polynomial that is linear in each variable. This work will give us a new perspective on the geometry of $X_M$.

Next, we will find a bound for the number of $\F_p$-rational points on configuration hypersurfaces (see Remark~\ref{rmk:other_bounds}), enhancing the existing bounds in the literature (see Remark~\ref{rmk:maye}). Readers primarily interested in the main result of this paper (Theorem~\ref{thm:estimate_of_limit}) can safely omit this section.

\begin{remark}\label{rmk:maye}
An estimate for $\#X_M(\F_p)$ can already be derived using results for complete intersections over finite fields found in~\cite{Lang-Weil, Hooley, Matera}. Since $X_M$ is geometrically irreducible when $M$ is irreducible, it meets the assumptions in ~\cite{Lang-Weil, Hooley, Matera}. For example, applying the Lang-Weil estimate~\cite{Lang-Weil} yields
\begin{equation*}
	\#X_M(\F_p)=p^{n-1}+O(p^{n-3/2}),
\end{equation*}
where $n=\#E-1$. To use the estimates from~\cite{Hooley, Matera}, we additionally need the dimension of the singular locus of $X_M$. It was shown in~\cite{Patterson, Denham-et21} that $X_M$ is highly singular. Consequently, the error term remains no better than the one given in this work (see Theorem~\ref{thm:rational_points_big_O}). 
\end{remark}

\subsection{Geometric aspects}

If $A$ is a commutative ring with unity and $f\in A$, then $A_f$ denotes the localization of $A$ at $f$. If in addition, $A$ is a graded ring and $f$ is a homogeneous element of $A$, then $A_{(f)}$ is the subring of elements of degree zero of $A_f$. We also fix the following notation.
\begin{itemize}
	\item $k$ : field.
	\item $\A^n_k\colonequals\spec k[T_1,\dots,T_n]$ : affine space of dimension $n$ over $k$.
	\item $\ps^n_k\colonequals\proj k[T_0,\dots,T_n]$ : projective space of dimension $n$ over $k$.
	\item $X_F\colonequals\proj k[T_0,\dots,T_n]/(F)$ : hypersurface cut out by a homogeneous polynomial $F$.
\end{itemize}
If, additionally, $\overline{G}\in k[T_0,\dots,T_n]/(F)$, then
\begin{itemize}
	\item $V_+(\overline{G})\colonequals\{\p\in\proj k[T_0,\dots,T_n]/(F)\colon \overline{G}\in \p\}$ : Zariski closed set of $X_F$, which we endow with the structure of closed subscheme of $X_F$ induced by the canonical morphism $\proj k[T_0,\dots,T_n]/(F,G)\hookrightarrow X_F$.
	\item $D_+(\overline{G})\colonequals X_F-V_+(\overline{G})$ : principal open set.
\end{itemize}
Finally, if $Y$ and $Z$ are closed subschemes of $\ps^n_k$, then $Y\cap Z$ is the scheme-theoretic intersection.

Let $F\in k[T_0,\dots,T_n]$ be a nonzero homogeneous polynomial that is linear in one of its variables; we can write (up to a permutation of the variables)
\begin{equation}\label{eq:1}
	F=T_0G_1+G_0,
\end{equation}
for some homogeneous polynomials $G_0,G_1\in k[T_1,\dots, T_n]$ with $G_1\neq 0$. Note that $\deg G_0=\deg F$ and $\deg G_1=\deg F-1$. 
Consider the hypersurfaces $X_F$, $X_{G_i}$, $i=0,1$. These come with canonical closed immersions $X_F\hookrightarrow\ps^n_k$,	 $X_{G_i}\hookrightarrow\ps^{n-1}_k$, $i=0,1$. Also, consider the subschemes $V_+(\overline{G}_1), D_+(\overline{G}_1)\subseteq X_F$.

\begin{proposition}\label{prop:iso_of_complement_of_G_1}
	$\ps^{n-1}_k-X_{G_1}$ and $D_+(\overline{G}_1)\subseteq X_F$ are isomorphic as $k$-schemes.
\end{proposition}
\begin{proof}
 	 By the universal property of localization, the  graded $k$-algebra homomorphism $ k[T_0,\dots,T_n]\to k[T_1,\dots, T_n]_{G_1}$, defined by $T_0\mapsto -\frac{G_0}{G_1}$ and $T_i\mapsto T_i$ for $i=1,\dots,n$, extends to a homomorphism of graded $k$-algebras $\varphi\colon k[T_0,\dots,T_n]_{G_1}\to k[T_1,\dots, T_n]_{G_1}$.

	We claim that $\ker\varphi=(F)_{G_1}$. The inclusion $(F)_{G_1}\subseteq \ker\varphi$ follows from observing that $\varphi(F)=(-G_0/G_1)G_1+G_0=0$. Suppose now $P/G_1^m\in \ker\varphi$ and write $P=T_0^lP_1+P_0$ for some $P_0,P_1\in k[T_1,\dots,T_n]$ and some nonnegative integer $l$. We have
	\begin{equation*}
		0=\varphi\left(\frac{P}{G_1^m}\right)=\frac{(-1)^lG_0^lP_1+G_1^lP_0}{G_1^{m+l}}\Rightarrow G_1^lP_0=(-1)^{l+1}G_0^lP_1,
	\end{equation*}
	from which it follows that
	\begin{align*}
		\frac{P}{G_1^m}&=\frac{T_0^lP_1+P_0}{G_1^m}=\frac{(T_0G_1)^lP_1}{G_1^{m+l}}+\frac{(-1)^{l+1}G_0^lP_1}{G_1^{m+l}}\\
		&=\left((T_0G_1)^l+(-1)^{l+1}G_0^l\right)\frac{P_1}{G_1^{m+l}}.
	\end{align*}
	Note that $(T_0G_1)^l+(-1)^{l+1}G_0^l=(F-G_0)^l+(-1)^{l+1}G_0^l$ and that the latter expression is divisible by $F$ (by the Binomial Theorem). Hence $P/G_1^m\in (F)_{G_1}$. This proves our claim. 
	
	It is clear that $\varphi$ is surjective, and it induces an isomorphism of graded $k$-algebras $\bar{\varphi}\colon \frac{k[T_0,\dots,T_n]_{G_1}}{(F)_{G_1}}\to k[T_1,\dots,T_n]_{G_1} $. Thus 
	\begin{equation*}
		\left( \frac{k[T_0,\dots,T_n]}{(F)}\right)_{\overline{G}_1}\cong \frac{k[T_0,\dots,T_n]_{G_1}}{(F)_{G_1}}\cong k[T_1,\dots,T_n]_{G_1},
	\end{equation*}
	where the first isomorphism is the canonical one. Since the degree zero components are preserved under these isomorphisms, we get
	\begin{equation}\label{eq:iso_algebras}
		\left( \frac{k[T_0,\dots,T_n]}{(F)}\right)_{(\overline{G}_1)}\cong k[T_1,\dots,T_n]_{(G_1)}.
	\end{equation}
	The rings in~\eqref{eq:iso_algebras} are the coordinate rings of the affine open sets $D_+(\overline{G}_1)$ and $\ps^{n-1}_k-X_{G_1}$, respectively; hence, the result follows.
\end{proof}

Let $p=[1,0,\dots,0]\in \ps^n_k$ be the rational point corresponding to the ideal $(T_1,\dots,T_n)$. Let $\pi_p\colon \ps^n_k-\{p\}\to \ps^{n-1}_k$ denote the projection centered at $p$ induced by the $k$-algebra homomorphism $\phi\colon k[T_1,\dots,T_n]\to k[T_0,\dots,T_n]$ given by $\phi(T_j)=T_j$ for all $j=1,\dots,n$. At the level of rational points we have $\pi_p([a_0,\dots,a_n])=[a_1,\dots,a_n]$.

\begin{corollary}\label{cor:projection_restricts_to_iso}
	The projection $\pi_p$ restricts to an isomorphism of $k$-schemes $D_+(\overline{G}_1)\to \ps^{n-1}_k-X_{G_1}$.
\end{corollary}
\begin{proof}
	Observe that $p\notin D_+(\overline{G}_1)$ and $\pi_p(D_+(\overline{G}_1))\subseteq \ps^{n-1}_k-X_{G_1}$. Thus, $\pi_p$ restricts to a morphism of $k$-schemes $ D_+(\overline{G}_1)\to \ps^{n-1}_k-X_{G_1}$. This morphism corresponds to a homomorphism $k[T_1,\dots,T_n]_{(G_1)}\to \left(\frac{k[T_0,\dots,T_n]}{(F)}\right)_{(\overline{G}_1)}$ that is induced by the localization of $\iota\circ \phi$ at $G_1$, where $\iota$ is the canonical homomorphism $k[T_0,\dots, T_n]\to k[T_0,\dots, T_n]/(F)$.
	
	It is straightforward to see that the localization of $\iota\circ \phi$ at $G_1$ and the homomorphism given in \eqref{eq:iso_algebras} are inverses of each other. Hence, the result follows.
\end{proof}

From Corollary~\ref{cor:projection_restricts_to_iso}, we obtain the following commutative diagram.
\begin{equation*}
	\begin{tikzcd}
		D_+(\overline{G}_1)\arrow[r, hook]\arrow[d, "\cong"'] &X_F-\{p\}\arrow[r, hook] & \ps^n_k-\{p\}\arrow[d, "\pi_p"]\\
		\ps^{n-1}_k-X_{G_1}\arrow[rr, hook]&&\ps^{n-1}_k
	\end{tikzcd}
\end{equation*}

When $F$ is irreducible, it follows from Proposition~\ref{prop:iso_of_complement_of_G_1} that $X_F$ is birational to $\ps^{n-1}_k$. Hence, $X_F$ is rational.

If $X=\proj k[T_1,\dots,T_n]/I$ is a projective scheme, then the scheme $\proj k[T_1,\dots,T_n][T_0]/I$ is called the  \emph{projective cone} of $X$ with \emph{vertex} $p$ and it is denoted by $C_p(X)$. It is well-known that $\pi_p$ induces a surjective morphism $\theta\colon C_p(X)-\{p\}\to X$ where the fibre of $\theta$ over $x\in X$ is isomorphic to $\A^1_{\kappa(x)}$ (here $\kappa(x)$ is the residue field of $x$); for more details we refer the reader to~\cite{Harris}*{Lecture $3$}.

\begin{proposition}\label{prop:intersection_of_G_1_and_G_zero}
	The projection morphism $\pi_p$ induces a surjective morphism of $k$-schemes $\theta \colon V_+(\overline{G}_1)-\{p\}\to X_{G_1}\cap X_{G_0}$, where the fibre of $\theta$ over $x\in X_{G_1}\cap X_{G_0}$ is isomorphic to $\A^1_{\kappa(x)}$.
\end{proposition}
\begin{proof}
	The closed immersion $X_F\hookrightarrow\ps^n_k$ induces an isomorphism of $k$-schemes $V_+(\overline{G}_1)\cong V_+(F)\cap V_+(G_1)$. In addition, $V_+(F)\cap V_+(G_1)=V_+(G_1,G_0)$ as schemes. Now notice that $G_1$ and $G_0$ are independent of $T_0$,  so that $V_+(G_1,G_0)=C_p(X_{G_1}\cap X_{G_0})$. Hence, the result follows.
\end{proof}

\subsection{Counting points over finite fields}

We continue using the notation from the previous subsection and assume $k$ is a finite field with $q$ elements.

\begin{proposition}\label{prop:rational_points}
	If $F\in k[T_0,\dots,T_n]$ is a homogeneous polynomial as in \eqref{eq:1}, then there is a natural identification 
	\begin{equation*}
		X_F(k)=C_p(X_{G_1}\cap X_{G_0})(k)\sqcup(\ps^{n-1}_k(k)-X_{G_1}(k)).
	\end{equation*}
\end{proposition}
\begin{proof}
	As topological spaces, $X_F=V_+(\overline{G}_1)\sqcup D_+(\overline{G}_1)$. In the proof of Proposition~\ref{prop:intersection_of_G_1_and_G_zero}, we showed that $V_+(\overline{G}_1)\cong C_p(X_{G_1}\cap X_{G_0})$. Furthermore, $D_+(\overline{G}_1)\cong \ps^{n-1}_k-X_{G_1}$ by Proposition~\ref{prop:iso_of_complement_of_G_1}. Both isomorphisms are over $k$; thus, the result follows.
\end{proof}

\begin{corollary}\label{cor:rational_points}
	If $F\in k[T_0,\dots,T_n]$ is a homogeneous polynomial as in \eqref{eq:1}, then
	\begin{equation*}
		\#X_F(k)=q\#(X_{G_1}(k)\cap X_{G_0}(k))+\#\ps^{n-1}_k(k)-\#X_{G_1}(k)+1.
	\end{equation*}
\end{corollary}
\begin{proof}
	It is a direct consequence of Proposition~\ref{prop:rational_points}.
\end{proof}

The equality in Corollary~\ref{cor:rational_points} was first deduced by Stembridge~\cite{Stembridge}. Nonetheless, the geometry behind this identity is concealed by the probabilistic methods used.

\begin{theorem}[\cite{Couvreur16}*{Corollary $3.3$}]\label{thm:couvreur}
	If $X\subseteq \ps^n_k$ is an equidimensional closed subscheme of dimension $d<n$ and degree $\delta$, then
	\begin{equation*}
		\#X(k)\leq \delta\Big(\#\ps^d_k(k)-\#\ps^{2d-n}_k(k)\Big)+\#\ps_k^{2d-n}(k).
	\end{equation*}
\end{theorem}
\noindent When $m$ is a negative integer then the value of $\#\ps^m_k(k)$ is zero by convention.

\begin{theorem}\label{thm:bounds_for_rational_points}
	Suppose that $F\in k[T_0,\dots,T_n]$ is a homogeneous polynomial as in \eqref{eq:1}. If $F$ is irreducible over $k$, then there are monic polynomials $f(t), g(t)\in\Z[t]$ of degree $n-1$ independent of $k$ such that
	\begin{equation*}
		g(q)\leq \#X_F(k)\leq f(q).
	\end{equation*}
	The coefficients of $g(t)$ and $f(t)$ depend only on $\deg F$ and $n$.
\end{theorem}

We need the following result on complete intersections to prove Theorem~\ref{thm:bounds_for_rational_points}. But before that, let us define what a complete intersection is. We say that the $k$-scheme $X=\proj k[T_0,\dots,T_n]/(f_1,\dots,f_r)$ is a \emph{complete intersection} if $\text{codim}(X,\ps^n_k)=r$ (equivalently, $\dim X=n-r$).

\begin{proposition}\label{prop:complete_intersection}
	Suppose that $F\in k[T_0,\dots,T_n]$ is a homogeneous polynomial as in \eqref{eq:1}. If $F$ is irreducible over $k$ of $\deg F>1$, then $X_{G_1}\cap X_{G_0}$ is  a complete intersection.
\end{proposition}
\begin{proof}
	We need to prove that $\text{codim}(X_{G_1}\cap X_{G_0},\ps^{n-1}_k)=2$. It is enough to show that the height of the ideal $(G_0,G_1)$, denoted $\text{ht}(G_0,G_1)$, is $2$ in $k[T_1,\dots,T_n]$. Since this ideal is generated by two polynomials, we have $\text{ht}(G_0,G_1)\leq 2$. Let $\p$ be a prime ideal containing $(G_0,G_1)$. Then $(G_0)\subseteq \p$. We claim that $\p$ is not a minimal prime ideal of $(G_0)$. Indeed, as $F$ is irreducible, $G_1$ and $G_0$ are coprime, therefore $\overline{G}_1$ is not a nonzero divisor of $k[T_0,\dots,T_n]/(G_0)$, which shows that $G_1$ cannot be contained in any minimal prime ideal of $G_0$. Consequently there must be a nonzero prime ideal $\mathfrak{q}$ strictly contained in $\p$. Then $\text{ht}(\p)\geq 2$, which implies that $\text{ht}(G_0,G_1)\geq 2$ as $\text{ht}(G_0,G_1)=\inf\{\text{ht}(\p)\colon (G_0,G_1)\subseteq \p\}$. This completes the proof.
\end{proof}

\begin{proof}[Proof of Theorem~\ref{thm:bounds_for_rational_points}]
	Let $m:=\deg F$. If $m=1$, then $X_F\cong \ps^{n-1}_k$ so that $\#X_F(k)=(q^n-1)/(q-1)$. Hence, we take $f(t)=g(t)=(t^n-1)/(t-1)$.
	
	Now assume $m>1$ so that $X_{G_1}\cap X_{G_0}$ is a complete intersection by Proposition~\ref{prop:complete_intersection}. Since complete intersections are equidimensional (see Section $10.135$ in \cite{stacks-project}), it follows that $X_{G_1}\cap X_{G_0}$ is equidimensional of dimension $n-3$. By B\'{e}zout's theorem (see Theorem III-$71$ in \cite{Eisenbud-harris00}) its degree is $m(m-1)$. We now use Theorem~\ref{thm:couvreur} to get
	\begin{align*}
		\#X_{G_1}(k)\cap X_{G_0}(k)\leq m(m-1)(\#\ps_k^{n-3}(k)-\#\ps^{n-6}_k(k))+\#\ps^{n-6}_k(k).
	\end{align*}
	The polynomial
	\begin{equation*}
		h(t)=m(m-1)\left(\frac{t^{n-2}-1}{t-1}\right)+(1+m-m^2)\left(\frac{t^{n-5}-1}{t-1}\right)\in\Z[t]
	\end{equation*}
	has degree $n-3$ and satisfies $\#(X_{G_1}(k)\cap X_{G_0}(k))\leq h(q)$. Moreover, the quantity $\#\ps^{n-1}_k(k)-\#X_{G_1}(k)+1$ is bounded above by $(q^n-1)/(q-1)+1$. Thus, if we take $f(t)=(t^n-1)/(t-1)+th(t)\in\Z[t]$, then $f(t)$ is monic of degree $n-1$ and $\#X_F(k)\leq f(q)$ by Corollary~\ref{cor:rational_points}.
	
	On the other hand, using Proposition~\ref{prop:rational_points} and applying Theorem~\ref{thm:couvreur} to $X_{G_1}$, we have
	\begin{align*}
		\#X_F(k)&\geq \#\ps^{n-1}_k(k)-\#X_{G_1}(k)+1\\
		&\geq \#\ps^{n-1}_k(k)-\Big((m-1)\Big(\#\ps^{n-2}_k(k)-\#\ps^{n-4}_k(k)\Big)+\#\ps^{n-4}_k(k)\Big)+1,
	\end{align*}
	so we define
	\begin{equation*}
		g(t)=\frac{t^n-1}{t-1}-(m-1)\left(\frac{t^{n-1}-1}{t-1} \right)+m\left( \frac{t^{n-3}-1}{t-1}\right)+1 \in\Z[t].
	\end{equation*}
	Then $g(t)$ is monic of degree $n-1$ and $g(q)\leq\#X_F(k)$.
\end{proof}

The final theorem of this subsection follows directly from Theorem~\ref{thm:bounds_for_rational_points}.
 
\begin{theorem}\label{thm:rational_points_big_O}
	Suppose that $F\in k[T_0,\dots,T_n]$ is a homogeneous polynomial satisfying \eqref{eq:1}. If $F$ is irreducible over $k$, then 
	\begin{equation*}
		\#X_F(k)=q^{n-1}+O(q^{n-2}).
	\end{equation*}
	The implied constant is computable and depends only on $\deg F$ and $n$.
\end{theorem}

\begin{remark}\label{rmk:other_bounds}
When $M$ is a nonempty irreducible regular matroid on $E$, $\Psi_M$ is irreducible over $\F_p$ (Proposition~\ref{prop:irreducibility_conf_poluynomial}) and satisfies the assumptions of Theorem~\ref{thm:rational_points_big_O} with $k=\F_p$. Hence, we have 
\begin{equation}\label{eq:estimate_X_M}
	\#X_M(\F_p)=p^{n-1}+O(p^{n-2}),
\end{equation}
where $n=\#E-1$. The implied constant is computable and depends only on $\#E$ and $r(M)$.

One possible approach to improving the error-term in this estimate involves utilizing the fact that  $X_{\Psi_{M\backslash e}}\cap X_{\Psi_{M/e}}$ is a complete intersection and examining its singular locus to apply the results in~\cite{Matera}. It is important to note that a bound for $X_{\Psi_{M\backslash e}}(\F_p)\cap X_{\Psi_{M/e}}(\F_p)$ was used to compute the error-term, but its singular locus was not taken into account.
\end{remark}

To conclude this section, we introduce the next proposition, which is the matroid version of Lemma $1.3$ and Corollary $1.4$ in~\cite{Aluffi-Marcolli09}. The proof given there works for this case too. We let $\Sigma_n$ denote the union of the coordinate hyperplanes in $\ps^{n-1}_k$, and $M^\ast$ denote the dual matroid of $M$.

\begin{proposition}\label{prop:aluffi-marcolli}
	If $M$ is a regular matroid on  $\#E=\{e_1,\dots,e_n\}$, then 
	\begin{equation*}
		\Psi_M(\lambda_{e_1},\dots,\lambda_{e_n})=\left( \prod_{i=1}^{n}\lambda_{e_i}\right)\Psi_{M^\ast}(\lambda_{e_1}^{-1},\dots,\lambda_{e_n}^{-1}). 
	\end{equation*}
	Moreover, the Cremona transformation induces an	 isomorphism of $k$-schemes $X_M-\Sigma_n\cong X_{M^\ast}-\Sigma_n$.
\end{proposition}		
			
			
\section{The Family $\{M_\lambda\}_{\lambda\in \N^{E(M)}}$}

In this section, we fix a prime number $p$. Let $M$ be a regular matroid on $E$. If $\lambda\in \N^E$, then we define $\lambda_p\in \N^E$ to be the map given by the rule $\lambda_p(e)=k$, where $1\leq k\leq p$ and $\lambda(e)\equiv k\Mod p$.

\begin{definition}
	The \emph{height} of a map $\lambda\in\N^E$ is $\height(\lambda)\colonequals\max\{\lambda(e) : e\in E\}$.
\end{definition}

\begin{definition}
	The \emph{density} of a subset $S\subseteq \N^E$ is 
	\begin{equation*}
		\mu(S)\colonequals\lim_{m\to\infty}\frac{\#\left(S\cap\{\lambda\in \N^E : \height(\lambda)\leq m\}\right)}{\#\{\lambda\in \N^E :  \height(\lambda)\leq m\}},
	\end{equation*}
	provided the limit exits.
\end{definition}

\begin{definition}
	We define $\J_p(M)\colonequals\{\lambda\in\N^E : p\mid\#\Jac(M_\lambda)\}$. If $G$ is a graph, we define $\J_p(G)\colonequals\J_p(M(G))$.
\end{definition}

When $M$ is a graphic matroid, i.e., $M=M(G)$ for some graph $G$, we see that $\J_p(M(G))=\{\lambda\in\N^{E(G)} \vcentcolon  p\mid\#\Jac(G_\lambda)\}$ by virtue of Proposition~\ref{prop:M_lambda_for_graphic_matroids}.

\begin{theorem}\label{thm:limit_of_height_of_lambdas}
	Let $M$ be a regular matroid on $E$. If $\#E=n$ and $\Psi_M\neq 1$, then
	\begin{equation}\label{eq:limit}
		\mu(\J_p(M))=\frac{(p-1)\#X_M(\F_p)+1}{(p-1)\#\ps_{\F_p}^{n-1}(\F_p)+1}.
	\end{equation}
\end{theorem}
\begin{proof}
	For $m\in\N$, we let
	\begin{align*}
		B_m&\colonequals\{\lambda\in\N^E \vcentcolon \height(\lambda)\leq m\}\\
		A_m&\colonequals\{\lambda\in B_m \vcentcolon \Psi_M(\lambda)\equiv0\Mod p\}.	
	\end{align*}
	By Theorem~\ref{thm:order_of_jacobian_parametrized_by_conf_polynomial}, $\Psi_M(\lambda)=\#\Jac(M_\lambda)$, therefore
	\begin{align*}
		\#A_m&=\#\{\lambda\in B_m : p\mid\#\Jac(M_\lambda)\}.
	\end{align*}
	Now suppose $m\geq p$ and write $m=pt_m+l$ for some $t_m\in\N$ and $0\leq l <p$. Consider the following map:
	\begin{align*}
		\theta_m\colon B_m &\to B_p\\
		\lambda &\mapsto \lambda_p.
	\end{align*}
	The map $\theta_m$ is surjective as $B_p\subseteq B_m$ and $\lambda_p=\lambda$ for all $\lambda\in B_p$. We will find a lower bound and an upper bound for the size of $\#\theta_m^{-1}(\gamma)$ with $\gamma\in B_p$. A map $\lambda\colon E\to \N$ is a preimage of $\gamma$ if and only if the following two conditions hold:
	\begin{enumerate}
		\item $1\leq\lambda(e)\leq pt_m+l$ for all $e\in E$ and
		
		\item for each $e\in E$, $\lambda(e)=pk_e+\gamma(e)$ for some nonnegative integer $k_e$.
	\end{enumerate}
	According to whether $\gamma(e)\leq l$ or $\gamma(e)> l$, the possible values for $k_e$ are:
	\begin{enumerate}
		\item[$(i)$] $0\leq k_e \leq t_m$, if $\gamma(e)\leq l$;
		
		\item[$(ii)$] $0\leq k_e \leq t_m-1$,  if $\gamma(e)> l$.
	\end{enumerate}
	It follows that
	\begin{equation*}
		t_m^n\leq \#\theta_m^{-1}(\gamma)\leq (t_m+1)^n.
	\end{equation*}
	Moreover, as $\theta_m^{-1}(A_p)=A_m$ we have
	\begin{align*}
		\#A_m=\sum_{\gamma\in A_p}\#\theta_m^{-1}(\gamma)
	\end{align*}
	then
	\begin{align*}
		t_m^n\#A_p\leq\#A_m\leq (t_m+1)^n\#A_p.
	\end{align*}
	On the other hand, $\#B_m=m^n=(pt_m+l)^n$ and $\#B_p=p^n$. Using the inequalities 
	\begin{align*}
		(pt_m+l)^n\leq (pt_m+p)^n=p^n(t_m+1)^n\quad\text{and}\quad
		p^nt_m^n\leq (pt_m+l)^n,
	\end{align*}
	we get 
	\begin{align*}
		\frac{\#A_m}{\#B_m}&\leq \frac{(t_m+1)^n\#A_p}{(pt_m+l)^n}\leq\frac{(t_m+1)^n\#A_p}{p^nt_m^n} =\frac{(t_m+1)^n}{t_m^n}\frac{\#A_p}{\#B_p},\\
		\frac{\#A_m}{\#B_m}&\geq\frac{t_m^n\#A_p}{(pt_m+l)^n}\geq \frac{t_m^n\#A_p}{(t_m+1)^np^n}=\frac{t_m^n}{(t_m+1)^n}\frac{\#A_p}{\#B_p}.
	\end{align*}
	So that,
	\begin{equation*}
		\left( \frac{t_m}{t_m+1}\right)^n \frac{\#A_p}{\#B_p}	\leq \frac{\#A_m}{\#B_m}\leq \left( \frac{t_m+1}{t_m}\right) ^n\frac{\#A_p}{\#B_p}.
	\end{equation*}
	By letting $m\to \infty$ on both sides and noting that $t_m\to\infty$, we get
	\begin{equation*}
		\mu(\J_p(M))=\lim_{m\to\infty}\frac{\#A_m}{\#B_m}=\frac{\#A_p}{\#B_p}.
	\end{equation*}
	Finally, notice that
	\begin{equation*}
		\frac{\#A_p}{\#B_p}=\frac{(p-1)\#X_M(\F_p)+1}{(p-1)\#\ps_{\F_p}^{n-1}(\F_p)+1}.
	\end{equation*}
	This follows by first noting that there is a bijection between $B_p$ and $\A^n_{\F_p}(\F_p)$. Next observe that there is a bijection between $A_p$ and $V(\Psi_M)(\F_p)\subseteq \A^n_{\F_p}(\F_p)$. Finally, note that each equivalence class of $\ps^{n-1}_{\F_p}(\F_p)$ contains $p-1$ points of $\A^n_{\F_p}(\F_p)$.
\end{proof}

\begin{example}\label{ex:cycle2}
	Let $G=\cycle_2$.
	\begin{figure}[H]
		\[\begin{tikzpicture}[scale=0.6,x=1.5cm, y=1.5cm,>=stealth',
			shorten > = 1pt,
			node distance = 3cm and 4cm,
			el/.style = {inner sep=2pt, align=left, sloped},
			every label/.append style = {font=\tiny}
			]
			\vertex[fill] (c1) at (180:1) {};
			\vertex[fill] (c2) at (0:1) {};
			\path 
			(c1) [bend left=45]edge node[el,above]  {$e$}   (c2)
			(c2) edge node[el,below] {$f$} (c1)
			;
		\end{tikzpicture}\]
		\caption{$\cycle_2$.}
	\end{figure}
	\noindent If $\lambda\in\N^{E(G)}$, then 	$G_\lambda=\cycle_{\lambda(e)+\lambda(f)}$, whence $\Jac(G_\lambda)\cong \Z/(\lambda(e)+\lambda(f))\Z$. Hence, $\{\Jac(G_\lambda)\}_{\lambda\in \N^{E(G)}}$ is the family of all finite cyclic groups.
	
	On the other hand, $\Psi_G=\lambda_e+\lambda_f$ and $X_G\cong \ps_{\F_p}^0$. Thus, by Theorem~\ref{thm:limit_of_height_of_lambdas} we obtain
	\begin{equation*}
		\mu(\J_p(G))=\frac{(p-1)\#X_G(\F_p)+1}{(p-1)\#\ps_{\F_p}^{1}(\F_p)+1}=\frac{1}{p},
	\end{equation*}
	as expected.
\end{example}

\begin{remark}
	The limit~\eqref{eq:limit} is completely determined if we know the quantity $\#X_M(\F_p)$. So, thinking of $\#X_M(\F_p)$ as a function of $p$, it is of interest to determine whether there is a polynomial relation for the values $\#X_M(\F_p)$ as $p$ varies. More specifically, consider the function $\#X_M\colon q\mapsto \#X_M(\F_q)$ defined on the set of prime powers $q$. One is interested in knowing whether $\#X_M\in\Z[q]$. In 1997, Kontsevich conjectured the following: for any graph $G$, $\#X_{M(G)}\in\Z[q]$, where $M(G)$ is the cycle matroid of $G$. In \cite{Stembridge}, Stembridge provided evidence in support of Kontsevich's conjecture by showing that the conjecture held for all graphs containing up to $12$ edges. Yet, Belkale and Brosnan \cite{Belkale} disproved the conjecture; in fact, they showed that these functions can be rather general. Loosely speaking, if $X$ is any arbitrary scheme of finite type over $\Z$, then the function $\#X: q\mapsto \#X(\F_q)$ can be written as  a suitable sum of functions $\#X_{M(G)}$, where $G$ is a graph. In this fashion, the best one can hope for is to be able to find bounds for $\#X_M(\F_q)$ independent of $q$ or describe the matroids for which $\#X_M\in\Z[q]$.
\end{remark}

\begin{theorem}\label{thm:estimate_of_limit}
	If $M$ is a nonempty irreducible regular matroid on $E$, then
	\begin{equation*}
		\mu(\J_p(M))=\frac{1}{p}+O\left( \frac{1}{p^2}\right).
	\end{equation*}
	The implied constant is computable and depends only on $r(M)$ and $\#E$.
\end{theorem}
\begin{proof}
	If $M$ is irreducible, then $\Psi_M$ is irreducible and it can be written as $\Psi_M=\lambda_e\Psi_{M\backslash e}+\Psi_{M/e}$ for any $e\in E$, where $\Psi_{M\backslash e}$ and $\Psi_{M/e}$ are independent of $\lambda_e$ (see Proposition~\ref{prop:deletion_contraction_conf_polynomial}). All the results from the previous section apply to $X_M$. In particular, 
	\begin{equation*}
		\#X_M(\F_p)=p^{n-1}+O(p^{n-2})
	\end{equation*}
	where $n=\#E-1$ (notice that $\#E\geq 3$ as $M$ is irreducible). On the other hand, Theorem~\ref{thm:limit_of_height_of_lambdas} says that
	\begin{equation*}
		\mu(\J_p(M))=\frac{(p-1)\#X_M(\F_p)+1}{(p-1)\#\ps_{\F_p}^n(\F_p)+1}.
	\end{equation*}
     By Theorem~\ref{thm:rational_points_big_O}, there exists a computable constant $C$ depending on $\#E$ and $r(M)$ such that
	\begin{equation}\label{bound:thm6.8}
		-Cp^{n-2}\leq\#X_M(\F_p)-p^{n-1}\leq Cp^{n-2}.
	\end{equation}
	 
	 To use the lower and upper bounds in (\ref{bound:thm6.8}) for $(p-1)\#X_M(\F_p)+1$, notice that
	\begin{align*}
		(p-1)\#X_M(\F_p)+1&=(p-1)\Big((\#X_M(\F_p)-p^{n-1})+p^{n-1}\Big)+1\\
		&=(p-1)(\#X_M(\F_p)-p^{n-1})+(p-1)p^{n-1}+1,
	\end{align*}
	therefore
	\begin{equation*}
		-(p-1)(Cp^{n-2}-p^{n-1})+1	\leq (p-1)\#X_M(\F_p)+1\leq (p-1)(Cp^{n-2}+p^{n-1})+1.
	\end{equation*}
	Dividing by $(p-1)\#\ps_{\F_p}^n(\F_p)+1=p^{n+1}$ yields
	\begin{equation*}
		\frac{1}{p}-\frac{C+1}{p^2}+\frac{C}{p^3}+\frac{1}{p^{n+1}}		\leq\frac{(p-1)\#X_M(\F_p)+1}{(p-1)\#\ps_{\F_p}^n(\F_p)+1}\leq	\frac{1}{p}+\frac{C-1}{p^2}-\frac{C}{p^3}+\frac{1}{p^{n+1}}.
	\end{equation*}
	Since $C\geq 1$, we finally get
	\begin{equation*}
		-\frac{C+1}{p^2}		\leq\frac{(p-1)\#X_M(\F_p)+1}{(p-1)\#\ps_{\F_p}^n(\F_p)+1}-\frac{1}{p}\leq	\frac{C+1}{p^2}.
	\end{equation*}
	This concludes the proof.
\end{proof}

The dependence of the constant from Theorem~\ref{thm:estimate_of_limit} on $\#E$ and $r(M)$ highlights that $\mu(\J_p(M))\approx 1/p$ is not what we should expect to hold for all values of $p$ and all nonempty irreducible regular matroids $M$. Let us consider the scenario where $\#E-r(M)$ is much larger than $p$. We would expect the probability that $\prod_{e\notin B}\lambda(e)\neq 0$ modulo $p$ to be $(1-p^{-1})^{\#E-r(M)}$, which would be close $0$. Therefore, $\mu(\J_p(M))$ should be very close to $1$. This particular situation is observed in graphs, as shown in the example below.

\begin{example}
	Let $G$ be the banana graph on $10$ edges. Its configuration polynomial is
	\begin{equation*}
		\Psi_G=\sum_{i=1}^{10}\prod_{\substack{j=1\\j\neq i}}^{10}T_j.
	\end{equation*}
It is not hard to show that for any prime $p$
\begin{equation*}
	\#X_G(\F_p)=p^9+35p^8-195p^7+510p^6-798p^5+798p^4-510p^3+195p^2-35p-1.
\end{equation*}
This allows us to use formula~(\ref{eq:limit}) and compute the precise value of $\mu(\J_p(G))$ for any prime $p$. In particular, we have that 
\begin{align*}
	\mu(\J_2(G))&=\frac{1013}{1024}\approx 0.9892\\
	&\\
	\mu(\J_3(G))&=\frac{53246}{59049}\approx 0.9017\\
	&\\
	\mu(\J_5(G))&=\frac{6305324}{9765625}\approx 0.6456.\\
\end{align*}

\end{example}

\begin{remark}
	Matroid theory comes with a dual theory, so that one can speak of the dual matroid of a  matroid $M$. If $M$ happens to be regular, then so is its dual. Furthermore, their Jacobians are isomorphic. 
	
	Starting with a regular matroid $M$, we obtain a new family of regular matroids $\{(M_\lambda)^\ast\}_{\lambda\in\N^{E(M)}}$, where $(M_\lambda)^\ast$ denotes the dual of $M_\lambda$. In this case, $\Jac(M_\lambda)\cong\Jac((M_\lambda)^\ast)$. Consequently,  Theorem~\ref{thm:limit_of_height_of_lambdas} and Theorem~\ref{thm:estimate_of_limit} can be used to estimate the proportion of Jacobians with nontrivial $p$-torsion in this new family. More concretely, if we define the set $\J_p^\ast(M)\colonequals\{\lambda\in\N^{E(M)} : p\mid \#\Jac((M_\lambda)^\ast)\}$, then
	\begin{equation*}
		\mu(\J_p^\ast(M))= \frac{(p-1)\#X_M(\F_p)+1}{(p-1)\#\ps_{\F_p}^{n-1}(\F_p)+1}.
	\end{equation*}
	If $M\neq \emptyset$ is irreducible, then
	\begin{equation*}
		\mu(\J_p^\ast(M))= \frac{1}{p}+O\left(\frac{1}{p^2} \right).
	\end{equation*}
\end{remark}

In general, given a regular matroid $M$ and a map $\lambda:E(M)\to\N$, the Jacobians of $M_\lambda$ and $M^\ast_\lambda$ are not isomorphic; in fact, their orders might not have any common prime factors. Also, it is not true that $\mu(\J_p(M))=\mu(\J_p(M^\ast))$.  Nonetheless,  Proposition~\ref{prop:aluffi-marcolli} allows us to relate the proportions of the Jacobians with nontrivial $p$-torsion in the families $\{M_\lambda\}_{\lambda\in\N^{E(M)}}$ and $\{M^\ast_\lambda\}_{\lambda\in\N^{E(M^\ast)}}$ as follows. 

Let us define $$\mathcal{S}_p(M)\colonequals\{\lambda\in\N^{E(M)} \vcentcolon p\mid \#\Jac(M_\lambda)~\text{and}~p\nmid \lambda(e)~\text{for all}~e\in E(M)\}.$$

Recall that $\Sigma_n$ denotes the union of the coordinate hyperplanes in $\ps^{n-1}_{\F_p}$. 

\begin{proposition}\label{prop:final}
	If $M$ is a regular matroid with $\#E(M)=n$, then
	\begin{equation*}
		\mu(\mathcal{S}_p(M))=\frac{(p-1)\#\Big(X_M(\F_p)\cap\big(\ps^{n-1}_{\F_p}(\F_p)-\Sigma_n(\F_p)\big)\Big)}{(p-1)\#\ps^{n-1}_{\F_p}(\F_p)+1}.
	\end{equation*}
\end{proposition}
\begin{proof}
	For $m\in\N$, we let
	\begin{align*}
		B_m&\colonequals\{\lambda\in\N^{E(M)}\colon \height(\lambda)\leq m\}\\
		D_m&\colonequals\{\lambda\in B_m\colon \Psi_M(\lambda)\equiv0\Mod p~\text{and}~p\nmid \lambda(e)~\text{for all}~e\in E(M)\}.	
	\end{align*}
	The argument of the proof of Theorem~\ref{thm:limit_of_height_of_lambdas} applies to $B_m$ and $D_m$ in place of $A_m$. Hence, one obtains
	\begin{equation*}
		\mu(\mathcal{S}_p(M))=	\lim_{m\to\infty}\frac{\#D_m}{\#B_m}=\frac{\#D_p}{\#B_p},
	\end{equation*}
	and
	\begin{equation*}
		\frac{\#D_p}{\#B_p}=\frac{(p-1)\#\Big(X_M(\F_p)\cap\big(\ps^{n-1}_{\F_p}(\F_p)-\Sigma_n(\F_p)\big)\Big)}{(p-1)\#\ps^{n-1}_{\F_p}(\F_p)+1}.
	\end{equation*}
\end{proof}

\begin{corollary}
	If $M$ is a regular matroid with $\#E(M)=n$, then $$\mu(\mathcal{S}_p(M))=\mu(\mathcal{S}_p(M^\ast)).$$
\end{corollary}
\begin{proof}
	By Proposition~\ref{prop:aluffi-marcolli}, $X_M-\Sigma_n\cong X_{M^\ast}-\Sigma_n$ over $\F_p$. Hence, the result follows from  Proposition~\ref{prop:final}.
\end{proof}		

\section*{Acknowledgements}		

This work was completed during my Ph.D. studies at the University of Western Ontario, and I wish to express my gratitude to the Mathematics Department for its warm hospitality and support.

Finally, I would like to thank the anonymous referees for their valuable insights, which greatly improved the quality of this paper.

\begin{bibdiv}
	\begin{biblist}
		\bib{Aluffi-Marcolli09}{article}{
	author={Aluffi, Paolo},
	author={Marcolli, Matilde},
	title={Feynman Motives of Banana Graphs},
	journal={Commun. Number Theory Phys.},
	volume={3},
	date={2009},
	number={1},
	pages={1--57},
	issn={1931-4523},
	review={\MR{2504753}},
	doi={10.4310/CNTP.2009.v3.n1.a1},
}
		
    \bib{Bacher-et97}{article}{
	author={Bacher, Roland},
	author={de la Harpe, Pierre},
	author={Nagnibeda, Tatiana},
	title={The Lattice of Integral Flows and the Lattice of Integral Cuts on
		a Finite Graph},
	language={English, with English and French summaries},
	journal={Bull. Soc. Math. France},
	volume={125},
	date={1997},
	number={2},
	pages={167--198},
	issn={0037-9484},
	review={\MR{1478029}},
}

\bib{Bake-Norine07}{article}{
	author={Baker, Matthew},
	author={Norine, Serguei},
	title={Riemann-Roch and Abel-Jacobi Theory on a Finite Graph},
	journal={Adv. Math.},
	volume={215},
	date={2007},
	number={2},
	pages={766--788},
	issn={0001-8708},
	review={\MR{2355607}},
	doi={10.1016/j.aim.2007.04.012},
}

\bib{Belkale}{article}{
	author={Belkale, Prakash},
	author={Brosnan, Patrick},
	title={Matroids, Motives, and a Conjecture of Kontsevich},
	journal={Duke Math. J.},
	volume={116},
	date={2003},
	number={1},
	pages={147--188},
	issn={0012-7094},
	review={\MR{1950482}},
	doi={10.1215/S0012-7094-03-11615-4},
}

\bib{Baker}{article}{
	author={Baker, Matthew},
	author={Norine, Serguei},
	title={Harmonic Morphisms and Hyperelliptic Graphs},
	journal={Int. Math. Res. Not. IMRN},
	date={2009},
	number={15},
	pages={2914--2955},
	issn={1073-7928},
	review={\MR{2525845}},
	doi={10.1093/imrn/rnp037},
}

\bib{Bloch}{article}{
	author={Bloch, Spencer},
	author={Esnault, H\'{e}l\`ene},
	author={Kreimer, Dirk},
	title={On motives associated to graph polynomials},
	journal={Comm. Math. Phys.},
	volume={267},
	date={2006},
	number={1},
	pages={181--225},
	issn={0010-3616},
	review={\MR{2238909}},
	doi={10.1007/s00220-006-0040-2},
}

\bib{Bosch02}{article}{
	author={Bosch, Siegfried},
	author={Lorenzini, Dino},
	title={Grothendieck's Pairing on Component Groups of Jacobians},
	journal={Invent. Math.},
	volume={148},
	date={2002},
	number={2},
	pages={353--396},
	issn={0020-9910},
	review={\MR{1906153}},
	doi={10.1007/s002220100195},
}

\bib{Kreimer}{article}{
	author={Broadhurst, D. J.},
	author={Kreimer, D.},
	title={Association of multiple zeta values with positive knots via
		Feynman diagrams up to $9$ loops},
	journal={Phys. Lett. B},
	volume={393},
	date={1997},
	number={3-4},
	pages={403--412},
	issn={0370-2693},
	review={\MR{1435933}},
	doi={10.1016/S0370-2693(96)01623-1},
}

\bib{Clancy-et15}{article}{
	author={Clancy, Julien},
	author={Kaplan, Nathan},
	author={Leake, Timothy},
	author={Payne, Sam},
	author={Wood, Melanie Matchett},
	title={On a Cohen-Lenstra Heuristic for Jacobians of Random Graphs},
	journal={J. Algebraic Combin.},
	volume={42},
	date={2015},
	number={3},
	pages={701--723},
	issn={0925-9899},
	review={\MR{3403177}},
	doi={10.1007/s10801-015-0598-x},
}

\bib{Couvreur16}{article}{
	author={Couvreur, Alain},
	title={An Upper Bound on the Number of Rational Points of Arbitrary
		projective varieties over finite fields},
	journal={Proc. Amer. Math. Soc.},
	volume={144},
	date={2016},
	number={9},
	pages={3671--3685},
	issn={0002-9939},
	review={\MR{3513530}},
	doi={10.1090/proc/13015},
}

\bib{Denham-et21}{article}{
	author={Denham, Graham},
	author={Schulze, Mathias},
	author={Walther, Uli},
	title={Matroid Connectivity and Singularities of Configuration
		Hypersurfaces},
	journal={Lett. Math. Phys.},
	volume={111},
	date={2021},
	number={1},
	pages={Paper No. 11, 67},
	issn={0377-9017},
	review={\MR{4205801}},
	doi={10.1007/s11005-020-01352-3},
}

\bib{Eisenbud-harris00}{book}{
	author={Eisenbud, David},
	author={Harris, Joe},
	title={The Geometry of Schemes},
	series={Graduate Texts in Mathematics},
	volume={197},
	publisher={Springer-Verlag, New York},
	date={2000},
	pages={x+294},
	isbn={0-387-98638-3},
	isbn={0-387-98637-5},
	review={\MR{1730819}},
}

\bib{Harris}{book}{
	author={Harris, Joe},
	title={Algebraic geometry},
	series={Graduate Texts in Mathematics},
	volume={133},
	note={A first course},
	publisher={Springer-Verlag, New York},
	date={1992},
	pages={xx+328},
	isbn={0-387-97716-3},
	review={\MR{1182558}},
	doi={10.1007/978-1-4757-2189-8},
}

\bib{Hooley}{article}{
	author={Hooley, C.},
	title={On the number of points on a complete intersection over a finite
		field},
	note={With an appendix by Nicholas M. Katz},
	journal={J. Number Theory},
	volume={38},
	date={1991},
	number={3},
	pages={338--358},
	issn={0022-314X},
	review={\MR{1114483}},
	doi={10.1016/0022-314X(91)90023-5},
}

\bib{Horton}{article}{
	author={Horton, Matthew D.},
	author={Stark, H. M.},
	author={Terras, Audrey A.},
	title={What are Zeta Functions of Graphs and What are They Good for?},
	conference={
		title={Quantum graphs and their applications},
	},
	book={
		series={Contemp. Math.},
		volume={415},
		publisher={Amer. Math. Soc., Providence, RI},
	},
	date={2006},
	pages={173--189},
	review={\MR{2277616}},
	doi={10.1090/conm/415/07868},
}

\bib{Lang-Weil}{article}{
	author={Lang, Serge},
	author={Weil, Andr\'e},
	title={Number of points of varieties in finite fields},
	journal={Amer. J. Math.},
	volume={76},
	date={1954},
	pages={819--827},
	issn={0002-9327},
	review={\MR{0065218}},
	doi={10.2307/2372655},
}

\bib{Lorenzini89}{article}{
	author={Lorenzini, Dino J.},
	title={Arithmetical Graphs},
	journal={Math. Ann.},
	volume={285},
	date={1989},
	number={3},
	pages={481--501},
	issn={0025-5831},
	review={\MR{1019714}},
	doi={10.1007/BF01455069},
}

\bib{Lorenzini08}{article}{
	author={Lorenzini, Dino},
	title={Smith Normal Form and Laplacians},
	journal={J. Combin. Theory Ser. B},
	volume={98},
	date={2008},
	number={6},
	pages={1271--1300},
	issn={0095-8956},
	review={\MR{2462319}},
	doi={10.1016/j.jctb.2008.02.002},
}

\bib{Marcolli}{book}{
	author={Marcolli, Matilde},
	title={Feynman motives},
	publisher={World Scientific Publishing Co. Pte. Ltd., Hackensack, NJ},
	date={2010},
	pages={xiv+220},
	isbn={978-981-4304-48-1},
	isbn={981-4304-48-4},
	review={\MR{2604634}},
}

\bib{Matera}{article}{
	author={Matera, Guillermo},
	author={P\'erez, Mariana},
	author={Privitelli, Melina},
	title={Explicit estimates for the number of rational points of singular
		complete intersections over a finite field},
	journal={J. Number Theory},
	volume={158},
	date={2016},
	pages={54--72},
	issn={0022-314X},
	review={\MR{3393540}},
	doi={10.1016/j.jnt.2015.06.007},
}

\bib{Merino}{thesis}{
	author={Merino, Criel},
	title={Matroids, the Tutte Polynomial and the Chip Firing Game},
	type={Ph.D. Thesis},
	school={University of Oxford},
	date={1999},
}
					
	\bib{Oxley11}{book}{
	author={Oxley, James},
	title={Matroid Theory},
	series={Oxford Graduate Texts in Mathematics},
	volume={21},
	edition={2},
	publisher={Oxford University Press, Oxford},
	date={2011},
	pages={xiv+684},
	isbn={978-0-19-960339-8},
	review={\MR{2849819}},
	doi={10.1093/acprof:oso/9780198566946.001.0001},
}

\bib{Patterson}{article}{
	author={Patterson, Eric},
	title={On the singular structure of graph hypersurfaces},
	journal={Commun. Number Theory Phys.},
	volume={4},
	date={2010},
	number={4},
	pages={659--708},
	issn={1931-4523},
	review={\MR{2793424}},
	doi={10.4310/CNTP.2010.v4.n4.a3},
}	

\bib{Stembridge}{article}{
	author={Stembridge, John R.},
	title={Counting Points on Varieties over Finite Fields Related to a
		Conjecture of Kontsevich},
	journal={Ann. Comb.},
	volume={2},
	date={1998},
	number={4},
	pages={365--385},
	issn={0218-0006},
	review={\MR{1774975}},
	doi={10.1007/BF01608531},
}	

\bib{stacks-project}{webpage}{
	author={The Stacks Project authors},
	title={The Stacks Project},
	url={https://stacks.math.columbia.edu},
	year={2021},
}	

\bib{Wood17}{article}{
	author={Wood, Melanie Matchett},
	title={The Distribution of Sandpile Groups of Random Graphs},
	journal={J. Amer. Math. Soc.},
	volume={30},
	date={2017},
	number={4},
	pages={915--958},
	issn={0894-0347},
	review={\MR{3671933}},
	doi={10.1090/jams/866},
}

\bib{Zapata}{thesis}{
	author={Zapata Ceballos, Sergio R},
	title={Distribution of the p-Torsion of Jacobian Groups of Regular Matroids},
	type={Ph.D. Thesis},
	school={Western University},
	date={2021},
}		
	\end{biblist}
\end{bibdiv}
\end{document}